\documentclass[12pt]{amsart}
\usepackage{amscd}
\usepackage{amssymb}
\usepackage{amsmath}
\usepackage{amsfonts}

\def\opn#1#2{\def#1{\operatorname{#2}}}
\opn\Spec{Spec}
\opn\deg{deg}
\opn\Coker{Coker}
\opn\Ker{Ker}
\opn\Pf{Pf}
\opn\ini{in}
\opn\Pf{Pf}

\def\mm{{\mathfrak m}}
\newtheorem{Theorem}{Theorem}[section]
\newtheorem{Proposition}[Theorem]{Proposition}
\newtheorem{Corollary}[Theorem]{Corollary}
\newtheorem{Lemma}[Theorem]{Lemma}
\newtheorem{Remark}[Theorem]{Remark}

\newcommand{\N}{\mathbb N}

\newcommand{\E}{\mathbb E}

\newcommand\alp{\alpha}

\newcommand\lam{\lambda}
\newcommand\var{\varphi}

\newcommand{\nor}{Noe\-the\-rian ring}

\newcommand{\rmo}{$R$-mo\-du\-le}

\newcommand{\wrt}{with respect to\ }

\newcommand{\fg}{fi\-ni\-te\-ly ge\-ne\-ra\-ted}

\newcommand{\CM}{Cohen-Macau\-lay}

\newcommand{\fgtfrmo}{fi\-ni\-te\-ly ge\-ne\-ra\-ted tor\-sion\-free
$R$-mo\-du\-le }

\newcommand{\bct}{\begin{center} \begin{tabular}}
\newcommand{\ect}{\end{tabular} \end{center}}
\newcommand{\pr}{^{\prime}}

\newcommand{\mat}{\begin{bmatrix}}
\newcommand{\emat}{\end{bmatrix}}
\newcommand{\edet}{\end{matrix} \right|}

\newcommand{\bca}{\begin{cases}}
\newcommand{\eca}{\end{cases}}

\renewcommand{\det}{\left| \begin{matrix}}

\newcommand{\col}{\colon}

\newcommand{\apl}[3]{#1 \col #2 \ra #3}
\newcommand{\hapl}[3]{#1 \col #2 \hra #3}

\newcommand{\oxr}{\otimes_{_{R}}}
\newcommand{\+}{\oplus}

\newcommand{\sube}{\subseteq}

\newcommand{\ra}{\rightarrow}

\newcommand{\hra}{\hookrightarrow}

\newcommand{\cP}{{\mathcal P}}
\newcommand{\cQ}{{\mathcal Q}}

\newcommand{\cS}{{\mathcal S}}

\newcommand{\cR}{{\mathcal R}}

\newcommand{\cF}{{\mathcal F}}

\newcommand{\cJ}{{\mathcal J}}

\newcommand{\frm}{{\mathfrak{m}}}

\newcommand{\frp}{{\mathfrak{p}}}
\newcommand{\frq}{{\mathfrak{q}}}

\DeclareMathOperator{\rank}{rank}

\DeclareMathOperator{\Ht}{ht}

\DeclareMathOperator{\Min}{Min}

\DeclareMathOperator{\Ass}{Ass}
\DeclareMathOperator{\ann}{ann}

\DeclareMathOperator{\z}{Z}

\DeclareMathOperator{\depth}{depth}
\DeclareMathOperator{\grade}{grade}

\DeclareMathOperator{\ad}{ad}

\DeclareMathOperator{\de}{d}

\DeclareMathOperator{\h}{H}

\newcommand{\sR}{\Spec(R)}

\newcommand{\dpt}{\depth \,}
\newcommand{\grd}{\grade \,}

\newcommand{\SE}{\cS(E)}

\newcommand{\RI}{\cR(I)}
\newcommand{\FI}{\cF(I)}

\newcommand{\RG}{\cR _G(G)}

\newcommand{\ol}{\overline}

\begin{document}

\title[Rees powers of a module]{Some asymptotic properties of \\ 
the Rees powers of a module}
\author{Ana L. Branco Correia}
\address{Centro de An\'{a}lise Matem\'{a}tica \\ Geometria e Sistemas 
Din\^{a}micos
\\ Instituto Superior T\'{e}cnico
\\ Av. Rovisco Pais \\ 1049-001 Lisboa
\\ Portugal
} \email{matalrbc@univ-ab.pt}
\address{ Universidade Aberta \\
Rua Fern\~ao Lopes $2^{\text{o}}$ Dto
\\ 1000-132 Lisboa \\ Portugal
}
\email{matalrbc@univ-ab.pt}
\thanks{}
\thanks{}
\author{Santiago Zarzuela}
\address{Departament d'\`Algebra i Geometria \\ Universitat de 
Barcelona \\ Gran Via 585 \\E-08007 Barcelona \\ Spain}
\email{zarzuela@mat.ub.es}
\thanks{Part of this work was made while the second author was visiting the
Instituto Superior T\'{e}cnico of the Universidade T\'{e}cnica de Lisboa,
to which he would like to thank for its hospitality.}
\thanks{The second author has been partially supported by BFM2001-3584
(Spain) and PR2003-0172 (Spain)} \subjclass{} \keywords{}
\date{}
\dedicatory{}
\commby{}
\begin{abstract}

Let $R$ be a commutative ring and let $G$ be a free $R$-module
with finite rank $e>0$. For any $R$-submodule $E\subset G$ one may
consider the image of the symmetric algebra of $E$ by the natural
map to the symmetric algebra of $G$, and then the graded
components $E_n$, $n\geq 0$, of the image, that we shall call the
$n$-th Rees powers of $E$ (with respect to the embedding $E\subset
G$). In this work we prove some asymptotic properties of the
$R$-modules $E_n$, $n\geq 0$, which extend well known similar ones
for the case of ideals, among them Burch's inequality for the
analytic spread.

\end{abstract}

\maketitle
\section{Introduction}
\noindent

Let $R$ be a commutative ring and let $G\simeq R^e$ be a free
$R$-module of finite rank $e>0$. For any $R$-submodule $E$ of $G$
we may consider the natural graded morphism from the symmetric
algebra of $E$ to the symmetric algebra of $G$ induced by the
embedding of $E$ in $G$, that we shall denote by $\phi :
\cS_R(E)\rightarrow \cS_R(G)$. As in the case of ideals, it is
natural in this situation to define the Rees algebra of $E$ (with
respect to the given embedding $E\subset G$) as the image of the
morphism $\phi$, namely $\cR _G(E):= \phi(\cS_R(E))\subset
\cS_R(G)$. Note that the symmetric algebra of $G$ is isomorphic to
the polynomial ring $R[t_1, \ldots , t_e]$ and that, by
definition, $\cR _G(G)=\cS_R(G)$. We may now consider the graded
components $[\cR _G(E)]_n$, for $n\geq 0$, of $\cR _G(E)$ that we
shall call the n-th Rees powers of $E$ (with respect to the given
embedding $E\subset G$) and denote by $E_n$, for $n\geq 0$. In
this way we may simply write $\cR _G(E)=\oplus _{n\geq 0} E_n$.
Observe that, by definition, we have that $E_n \subset G_n \simeq
R^g$, with $g={\binom{n+e-1}{e-1}}$ for $n\geq 0$.

\smallskip

Some general properties of the modules $E_n$, extending similar
ones for ideals, have been proved for instance by D. Katz and C.
Naude \cite{KaNa} and in the two dimensional case by D. Katz and
V. Kodiyalam \cite{KaKo}, V. Kodiyalam \cite{Ko} and R. Mohan
\cite{Mo}. In particular, it is shown in \cite{KaNa} that if $R$
is a Noetherian ring then the sets of associated primes $\Ass
(G_n/E_n)$ stabilize for $n \gg 0$, that is, $\Ass (G_n/E_n) =
\Ass (G_{n+1}/E_{n+1})$ for $n \gg 0$, the corresponding result
for ideals having been proved by M. Brodmann \cite{Br1}.

\smallskip

One of the most basic asymptotic properties of the powers of an
ideal is given by the well known Burch's inequality relating the
analytic spread of an ideal $I\subset R$ with the depths of the
$R$-modules $R/I^n$, for $n\geq 0$. Recall that if $(R, \mm)$ is a
local ring of dimension $d$, then the analytic spread $\ell(I)$ of
$I$ may be defined as the dimension of the ring $\cF (I)=\cR
(I)/\mm \cR (I)$, the fiber cone (or special fiber) of $I$. It is
well known that $\Ht (I)\leq \ell(I) \leq d$, and L. Burch
\cite{Bu} proved that $\ell(I)\leq d - \inf _{n\geq 1} \depth
R/I^n$. It is also known that this inequality becomes equality
when the associated graded ring of $I$ is Cohen-Macaulay, see for
instance D. Eisenbud and C. Huneke \cite{EiHu}.

\smallskip

Our main result in this paper is a natural extension to the case
of modules of the Burch's inequality. Namely, we shall be able to
prove the following result.

\begin{Theorem}
Let $(R, \mm )$ be a local ring and let $G\simeq R^e$ be a free
$R$-module of finite rank $e>0$. Assume that $\dim R=d>0$. Let
$E\subset G$ be an $R$-submodule of $G$ such that $E\neq G$ and
denote by $\ell _G(E)$ the dimension of the graded ring $\cR
_G(E)/\mm \cR _G(E)$. Then
$$\ell _G(E)\leq d + e - 1 - \inf _{n\geq 1} \depth G_n/E_n.$$
\end{Theorem}

We may also prove that the above inequality becomes equality
if the ring $R$ and the Rees algebra $\cR _G(E)$ are both
Cohen-Macaulay. In the case of modules, there is a lack of a
similar object as the associated graded ring of an ideal, but it
is well known that if $R$ is Cohen-Macaulay and $I\subset R$ is an
ideal whose Rees algebra $\cR (I)$ is Cohen-Macaulay, then the
associated graded ring of $I$ is Cohen-Macaulay too.

\smallskip

In order to prove Burch's inequality for modules we follow the
approach by M. Brodmann \cite{Br2} in the ideal case, which is
based on the stable behavior of the depths of $R/I^n$, for $n >
0$. Roughly speaking, by using the result proven by D. Katz and C.
Naude on the stability of the sets of associated primes of the
$R$-modules $G_n/E_n$, for $n \gg 0$, we can first extend to the
case of modules the result by M. Brodmann on the stable behavior
of the depths of $G_n/E_n$, for $n \gg 0$, and then to prove
Burch's inequality for modules.

\smallskip

Although in the case of ideals one can quickly prove Burch's
inequality by using the associated graded ring, we are forced to
follow this more complicated approach due to the lack of such a
similar object in the theory of Rees algebras of modules.

\smallskip

As a final application we also extend to modules a famous criteria
proved by R. C. Cowsik and M. V. Nori \cite{CoNo} for an ideal to
be a complete intersection. A special case of complete
intersection modules was already considered by D. A. Buchsbaum and
D. Rim in \cite{BuRi} under the name of parameter matrices,
playing the same role as system of parameters in their theory of
multiplicities for modules of finite length (nowadays called
Buchsbaum-Rim multiplicities). Other authors like D. Katz and C.
Naude \cite{KaNa} have also studied the properties of complete
intersection modules (or modules of the principal class).

\smallskip

Throughout this paper we shall always assume that $(R, \mm)$ is a
local ring of dimension $d$ with maximal ideal $\mm$ and that
$G\simeq R^e$ is a free $R$-module with finite rank $e>0$. For a
given $R$-submodule $E\subset G$ of $G$ we shall understand that
an embedding of $E$ in $G$ has been fixed. Consequently, the Rees
algebra $R_G(E)$ of $E$, and so the $n$-th Rees powers of $E$,
will always be with respect to this fixed embedding. One should
note that, for a given $R$-module $E$ and different embeddings of
$E$ into free $R$-modules, we can get non isomorphic Rees algebras
of $E$, see for instance A. Micali \cite[Chapitre III, 2. Un
example]{Mi} or the more recent D. Eisenbud, C. Huneke and B.
Ulrich \cite[Example 1.1]{EiHuUl}. See also these papers for a
discussion about the uniqueness of the definition of the Rees
algebra of a module. In particular, it is known that if an
$R$-module $E$ has rank, that is, $E\otimes _R Q$ is free with
positive rank, where $Q$ is the total ring of fractions of $R$,
then for any embedding $E\subset G\simeq R^e$ of $E$ into a free
$R$-module of positive rank, the Rees algebra of $E$ (with respect
to this embedding) is isomorphic to $\cS_R(E)/T_R(\cS_R(E))$, the
symmetric algebra of $E$ modulo its $R$-torsion, see
\cite{EiHuUl}.

\smallskip

We close this introduction with the following observation. Given
$E \subset G \simeq R^e$ an \rmo{} with $E \neq G$, the ideal of
$\RG$ generated by $E$ is the graded ideal
$$ E \RG = E \+ E \cdot G \+ E \cdot G_2 \+ \cdots =
\bigoplus_{i \geq 0} E \cdot G_i.$$ For each $n \geq 1$,
\begin{equation} \label{feR}
    (E \RG)^{n} =\bigoplus_{i \geq 0} E_{n} \cdot G_i
= E_{n} \RG
\end{equation}
and so, in particular, $ [E_n \RG]_{n} = E_n$.

\section{The analytic spread of a module}

Let $E\subset G$ be an $R$-submodule of $G$. If $E$ has rank $e>0$
there is an explicit formula for the dimension of
the Rees algebra of $E$, namely $\dim \cR _G(E)= d + e$,
see for instance A.
Simis, B. Ulrich and W. V. Vasconcelos \cite[Proposition
2.2]{SiUlVa}. In the general case, there is not such an explicit
formula, but one can prove that $\dim \cR _G(E)\leq d+e$. For this
we follow similar steps as in the proof of \cite[Proposition
2.2]{SiUlVa}. First, we determine the set of minimal primes of the
Rees algebra of $E$.

\begin{Lemma}\label{mpG}
    Let $E \subset G \simeq R^e$, $e>0$,
    be an \rmo. Then
    $$ \Min \cR_G(E) = \{ \cP =\frp \RG \cap \cR_{G}(E) \mid
    \frp \in \Min R \}. $$
\end{Lemma}

\begin{proof}
    For any $R$-ideal $J$, we set $\cJ= J \RG \cap \cR_{G}(E)$. It
    is easy to prove that if $(0) = \frq_{1} \cap \ldots \cap
    \frq_{s}$ is a shortest primary decomposition in $R$ then $(0) =
    \cQ_{1} \cap \ldots \cap \cQ_{s}$ is a shortest primary decomposition
    in $\cR_{G}(E)$. Hence
    $$ \Ass \cR_{G}(E) = \{ \sqrt{\cQ_{1}},\ldots,\sqrt{\cQ_{s}}
    \} \supseteq \Min \cR_{G}(E). $$
    Now, if $\cP \in \Min \cR_{G}(E)$, $\cP= \sqrt{\cQ_{i}}$ for
    some $1 \leq i \leq s$. But $$\sqrt{\cQ_{i}} = \sqrt{\frq_{i}
    \RG \cap \cR_{G}(E)} = \sqrt{\frq_{i}} \RG \cap \cR_{G}(E)$$
    (since $\RG$ is a polynomial ring) and, by minimality of
    $\sqrt{\cQ_{i}}$, $\sqrt{\frq_{i}} \in \Min (\Ass R)= \Min R$.
    For the other inclusion, let $\cP= \frp \RG \cap \cR_{G}(E)$ with
    $\frp \in \Min R$.
    Then $\frp = \sqrt{\frq_{i}}$ for some $1 \leq i \leq s$
    and we have $\cP = \sqrt{\cQ_{i}}= \sqrt{\frq_{i}
    \RG \cap \cR_{G}(E)} \in \Min \cR_{G}(E)$. The equality
    follows.
\end{proof}

\begin{Lemma}\label{mpG1}
    Let $E \subset G \simeq R^e$, $e>0$,
    be an \rmo. Let $\frp \in \sR$ and set $\cP = \frp \RG \cap
    \cR_{G}(E)$. Then
    $$ \cR_{G}(E)/\cP \simeq \cR_{\ol{G}}(\ol{E}), $$
    where $\ol{R} = R/\frp$, $\ol{G}=G/\frp G \simeq G \oxr \ol{R}$,
    $\ol{E}= (E+\frp G)/\frp G \subset \ol{G}$.
\end{Lemma}

\begin{proof}
    Let $\apl{\pi}{G}{\ol{G}}$ be the canonical epimorphism. Since
    $\pi(E) = (E+\frp G)/\frp G= \ol{E}$, the map $\apl{\pi\pr}
    {E}{\ol{E}}$ defined by $\pi\pr(z)=\pi(z)$ for every $z \in
    E$ is an epimorphism. Hence
    $\apl{\cS(\pi\pr)}{\cS_R(E)}{\cS_{\ol{R}}(\ol{E})}$
    is also an epimorphism. Therefore there exists an
    $R$-epimorphism $\apl{\rho}{\cR_G(E)}{\cR_{\ol{G}}(\ol{E})}$.
    Moreover, if $\hapl{\iota}{\cR_{\ol{G}}(\ol{E})}{\cR_{\ol{G}}(\ol{G})}$
    denotes the natural inclusion and $\apl{\lam}{\RG}{\RG/\frp \RG
    \simeq \cS_{\ol{R}}(\ol{G}) = \cR_{\ol{G}}(\ol{G})}$ denotes the canonical
    epimorphism, then
    $$ \lam|_{\cR_G(E)} = \iota \circ \rho.$$
    It follows that
    $$ \cP = \frp \RG \cap \cR_G(E) = \ker \lam|_{\cR_G(E)}
    = \ker (\iota \circ \rho) = \rho^{-1}(\ker \iota) = \ker \rho$$
    and so
    $$ \cR_G(E)/\cP \simeq \cR_{\ol{G}}(\ol{E})$$
    as claimed.
\end{proof}

Our bound for the analytic spread of a module will be a
consequence of the following expression of the dimension of the
Rees algebra.

\begin{Proposition}\label{mpG2}
    Let $E \subset G \simeq R^e$, $e>0$,
    be an \rmo. Then
    $$ \dim \cR_G(E) = \max \{ \dim R/\frp + \rank (E+\frp
    G/\frp G) \mid \frp \in \Min R \} \leq d + e. $$
\end{Proposition}

\begin{proof}
    We have, by the previous lemmas,
    \begin{align*}
     \dim \cR_G(E) & = \max \{ \dim \cR_G(E)/\cP \mid \cP \in
     \Min \cR_G(E) \} \\
     & = \max \{ \dim \cR_{G/\frp G}(E+\frp G/\frp G) \mid \frp \in
     \Min R \}.
    \end{align*}
    On the other hand, for each $\frp \in \Min R$, $\ol{R}=R/\frp$ is
    a domain and so $\ol{E}= E+\frp G/\frp G \sube G/\frp G =
    \ol{G} \simeq \ol{R}^e$ is a \fgtfrmo{} having rank $r \leq e$.
    Hence, in this case we have
    $$\dim \cR_{G/\frp G}(E+\frp G/\frp G)= \dim \ol{R} +
    \rank \ol{E} \leq d + e $$
    and the result follows.
\end{proof}

\begin{Remark}
{\rm Observe that the above result and its proof are also valid
for any Noetherian ring $R$, not necessarily local.}
\end{Remark}

As in the ideal case, we may define the fiber cone (or special
fiber) of $E$ as $\cF _G(E):= \cR _G(E)/\mm \cR_G (E)$, and the
analytic spread of $E$ as $\ell _G(E):= \dim \cF _G(E)$. From the
above proposition we get the following bound for the analytic
spread.

\begin{Corollary}\label{mpG3}
    Let $E \subset G \simeq R^e$, $e>0$,  be
    an \rmo. Assume $d>0$. Then
    $$ \ell_G(E) \leq d+e-1. $$
\end{Corollary}

\begin{proof}
    By definition
    $$ \ell_G(E) = \dim \cF_G(E) = \dim \cR_G(E)/
    \frm \cR_G(E). $$
    If $\frm \cR_G(E) \sube \cP$ for some $\cP \in \Min
    \cR_G(E)$, then
    $$ \frm = \frm \cR_G(E) \cap R \sube \cP \cap R = \frp \RG
    \cap \cR_G(E) \cap R = \frp \RG \cap R = \frp$$
    where $\frp \in \Min R$ (by Lemma \ref{mpG}). But $\dim R >0$,
    that is $\frm \not \in \Min R$, and so  $\frm \cR_G(E)$ is not
    contained in any minimal prime of $\cR_G(E)$. It follows that
    $\Ht \frm \cR_G(E) > 0$ and so
    $$ \ell_G(E) = \dim \cR_G(E)/\frm \cR_G(E) <
    \dim \cR_G(E) \leq d+e,$$
    proving the asserted inequality.
\end{proof}

Although in this paper we are not going to use the theory of
reductions of modules, it is worthwhile to point out that, under
suitable conditions similar to the ideal case, the analytic spread
of an $R$-module $E\subset G$ coincides with the minimal number of
generators of any minimal reduction of $E$. In particular, and if
$E$ has rank, it also holds that $\rank E\leq \ell _G(E)$, see
\cite[Propositon 2.2]{SiUlVa}.

\section{The asymptotic behavior of $\depth G_n/E_n$}

Our aim in this section is to prove that $\depth G_n/E_n$ takes a
constant value for large $n$. For that we shall need the following
technical lemma.

\begin{Lemma}\label{ad1}
    Let $E \varsubsetneq G \simeq R^e$, $e>0$, be
    an \rmo{} and let $a \in \frm$.
    Then, for every $n \geq 1$,

    \begin{enumerate}
    \item  $\ol{E}_n \simeq (E_n + a G_n) / a G_n$;

    \item  $ (G_n/E_n)/ a (G_n/E_n) \simeq \ol{G}_n /
    \ol{E}_n$;

    \item  $[\cF_{\ol{G}}(\ol{E})]_{n}
    \simeq E_n/ \frm E_n + (aG_n \cap E_n)$;

    \end{enumerate}
    where $\ol{E}=E+aG/aG\subset \ol{G}=G/aG$.
\end{Lemma}

\begin{proof}
    Clearly, $\ol{E} \neq \ol{G}$
    (by Nakayama's Lemma) since $E \neq G$ and $a \in \frm$.
    Moreover, for each $n \geq 1$, we have the commutative diagram

    $$ \begin{array}{ccccl}
        \cS(\ol{E})_{n} & \longrightarrow &
        \cS(\ol{G})_{n} & = & \ol{G}_n  \\
        \simeq &  & \simeq &&\\
        \SE_{n} \oxr R/(a) & \longrightarrow & \cS(G)_{n} \oxr R/(a) &
        = & G_n \oxr R/(a) \\
        \simeq &  & \simeq && \\
        \SE_{n}/a \SE_{n} & \longrightarrow & \cS(G)_{n}/a \cS(G)_{n} & = & G_n /a G_n
        \end{array}.  $$
    Hence, by the definition of the $n$-th Rees powers, it follows that
    $$ \ol{E}_n \simeq (E_n + a G_n) / a G_n, $$
    and (1) is proved.
    Moreover,
    \begin{align*}
        (G_n/E_n)/ a (G_n/E_n)  & \simeq G_n / (E_n + a G_n) \\
    & \simeq (G_n/a G_n)/ ((E_n + a G_n) / a G_n) \\
    & \simeq \ol{G}_n / \ol{E}_n
    \end{align*}
    and (2) is proved. As for (3) we have
    \begin{align*}
        [\cF_{\ol{G}}(\ol{E})]_{n} & =
    \ol{E}_n / \frm \ol{E}_n \simeq ((E_n + a G_n) / a G_n)/\frm (
    (E_n + a G_n) / a G_n) \\
    & \simeq (E_n + a G_n)/ (\frm E_n + a G_n)
    \simeq E_n / ((\frm E_n + a G_n) \cap E_n) \\
    & \simeq E_n / (\frm E_n +
    (a G_n \cap E_n))
    \end{align*}
    - the last isomorphism by the modular law.
\end{proof}

\begin{Remark}\label{an1}
   {\rm The lemma above still holds if $a \in \frm$ is replaced by an
    $R$-ideal $I \sube \frm$.}
\end{Remark}

For each $n \geq 1$, we denote the associated prime ideals of
$G_n/E_n$ by $A(n)$. By \cite[Theorem~2.1]{KaNa}, $A(n)$ is stable
for large $n$.

\begin{Theorem}\label{ad2}
    Let $E \varsubsetneq G \simeq R^e$, $e>0$,
    be an \rmo{}.
    Then $\dpt G_n/E_n$
    takes a constant value for large $n$.
\end{Theorem}

\begin{proof}
    We use induction on the inferior limit
    (\footnote{The inferior limite is by definition the smallest of the
    sublimits.})
    $$ \alp = \liminf_{n \ra \infty} \dpt G_n/E_n. $$
    If $\alp = 0$ then there exits $m \gg 0$ such that $\dpt G_m/
    E_m = 0$. Hence $\frm \in \Ass G_m/E_m = A(m)$.
    But $A(n)$ is stable by
    large $n$. Hence $\frm \in A(n)$ for all $n
    \geq m$, and so $\dpt G_n/E_n = 0$ for all $n \geq m$.

    Now suppose that $\alp >0$. Hence there exists $m \geq 1$ such
    that $\dpt G_n /E_n \neq 0$ for all $n \geq m$, that is
    $\frm \not \in A(n)$ for all $n \geq m$.
    Therefore, there
    exists $a \in \frm$ such that $a \not \in \bigcup_{\frp \in A(n)}
    \frp$ for sufficiently large $n$.
    It follows that
    $$ a \not \in \z_{R}(G_n/E_n) = \displaystyle
    {\bigcup_{\frp \in A(n)} \frp} $$
    and we have
    $$ \dpt (G_n/E_n)/a(G_n/E_n) = \dpt G_n/E_n
    -1 \; \; (n \gg 0).$$
    On the other hand, by the previous lemma,
    $$ (G_n/E_n)/ a (G_n/E_n) \simeq_{R/(a)} \ol{G}_n /
    \ol{E}_n $$
    and so
    $$ \beta = \liminf_{n \ra \infty} \dpt \ol{G}_n /
    \ol{E}_n < \liminf_{n \ra \infty}
    \dpt G_n/E_n = \alp. $$
    By induction hypothesis, $\dpt \ol{G}_n / \ol{E}_n$ takes a
    constant value for $n \gg 0$. Hence the result follows by
    induction.
\end{proof}

In the following, we shall denote by $\dpt (G, E)$ this asymptotic
constant value of $\dpt G_n/E_n$.

\section{The asymptotic behavior of the analytic spread of a
module}

In this section we prove our result extending to modules the
Burch's inequality, see \cite{Bu}. To do this, we need the
following lemmas. \smallbreak

\begin{Lemma}\label{bi1}
    Let $E \varsubsetneq G
   \simeq R^e$, $e>0$, be an \rmo{}, and assume $d>0$.
    Then, for all $n$,
    $$ \bigoplus_{i \geq 1} \left[
    \sqrt{\ann_{\cR_G(E)}(\cF_G(E))} \, \right]_{i}=
    \displaystyle{\bigoplus_{i \geq 1} \left[ \sqrt{\ann_{\cR_G(E)}
    (\+_{m \geq n} E_{m} /\frm E_{m})} \, \right]_{i}}. $$
\end{Lemma}

\begin{proof}
    The inclusion $\sube$ is clear. For the other one it suffices
    to see that for any homogeneous element $a \in E_{h}$ $(h>0)$
    such that $a E_{m}\sube \frm E_{hm}$ for any $m \geq n$, there exists
    $s>0$ with $a^s\in \frm R_G(E)$. Let $s>0$ such that $(s-1)h\geq n$.
    Then, $a^s=aa^{s-1}\sube aE_{(s-1)h}\sube \frm E_{sh}$. Hence
    $a^s\in \frm R_G(E)$ and the lemma follows.
    \end{proof}

\begin{Remark}\label{bi2}
    {\rm Let $E \varsubsetneq G \simeq R^e$, $e>0$, be an \rmo{}.
    Suppose that $d>0$ and $\dpt (G,E)>0$. Then there exists $a \in \frm$ such
    that $a \not \in \z_{R}(G_n/E_{n})$ for all $n \gg 0$ and $
    a \not \in \frp$ for all $\frp \in \Min R$.}
\end{Remark}

\begin{proof}
    Since $\dpt (G,E)>0$, then $\alp = \liminf_{n \ra \infty} \dpt
    G_n/E_{n} >0$ and, as in the proof of Theorem \ref{ad2},
    $\frm \not \in A(n)$ for all $n \gg 0$. Moreover, $\frm \not \in
    \Min R$, since $\dim R >0$. Hence
    $$ \frm \not \in \left(\bigcup_{n \gg 0} A(n)\right) \cup \Min R, $$
    It follows that there exists $a \in \frm$ such
    that $a \not \in \z_{R}(G_n/E_{n})$ for all $n \gg 0$ and
    $a \not \in \frp$ for all $\frp \in \Min R$, as asserted.
\end{proof}

\begin{Lemma}\label{bi3}
    Let $E \varsubsetneq G \simeq R^e$, $e>0$, be an \rmo{}.
    Suppose that $d>0$ and $\dpt (G,E)>0$. Then
    $$ \ell_G(E) = \ell_{\ol{G}}(\ol{E}), $$
    where $\ol{E}=E + aG /aG$, $\ol{G}=G/aG$ and
    $a \in \frm \setminus \bigcup_{n \gg 0}\z_{R}(G_n/E_{n})$.
\end{Lemma}

\begin{proof}
    Since $\dpt (G,E)>0$, we may choose an
    $a \in \frm \setminus \bigcup_{n \gg 0}\z_{R}(G_n/E_n)$.
    By Lemma \ref{ad1}, we have, for each $n$,
    $$[\cF_{\ol{G}}(\ol{E})]_{n}  \simeq E_n/
    (\frm E_n + (aG_n \cap E_n)). $$
    Now, since $a \in \frm$ and $a$ is regular \wrt $G_n/E_n$
    for all sufficiently large $n$
    $$(E_n :_{G_n}^{\ } a) = \{ z \in G_n
    \mid az \in E_n \} = E_n$$
    for all $n \gg 0$. Hence, for $n \gg 0$,
    $a G_n \cap E_n = a E_n \sube \frm E_n$,
    and so $\frm E_n + (aG_n \cap E_n) = \frm E_n$.
    It follows that
    $$ [\cF_{\ol{G}}(\ol{E})]_{n} \simeq [\cF_{G}(E)]_{n} $$
    for $n \gg 0$.
    Therefore, by the Lemma \ref{bi1}, we have for $n \gg 0$
    \begin{align*}
        \bigoplus_{i \geq 1} \left[\sqrt{\ann_{\cR_{G}(E)}
    \cF_{G}(E)} \, \right]_{i}
    & = \bigoplus_{i \geq 1} \left[\sqrt{\ann_{\cR_{G}(E)}
    ( \+_{m \geq n} \cF_{G}(E)_{m})} \, \right]_{i} \\
    & = \bigoplus_{i \geq 1} \left[\sqrt{\ann_{\cR_{G}(E)}
    ( \+_{m \geq n} [\cF_{\ol{G}}(\ol{E})]_{m})} \, \right]_{i} \\
    & = \bigoplus_{i \geq 1} \left[
    \sqrt{\ann_{\cR_{\ol{G}}(\ol{E})}
    ( \+_{m \geq n} [\cF_{\ol{G}}(\ol{E})]_{m})} \, \right]_{i} \\
    & = \bigoplus_{i \geq 1} \left[
    \sqrt{\ann_{\cR_{\ol{G}}(\ol{E})}
    \cF_{\ol{G}}(\ol{E})} \, \right]_{i}.
    \end{align*}
    Moreover, since $a \in \frm$,
    \begin{align*}
        [\ann_{\cR_{G}(E)} \cF_{G}(E)]_{0} & = \ann_{R} R/\frm =
    \ann_{R/(a)} (R/(a))/ \frm (R/(a)) \\
    & = [\ann_{\cR_{\ol{G}}(\ol{E})}
    \cF_{\ol{G}}(\ol{E})]_{0}.
    \end{align*}
    Now, since
    \begin{align*}
        \ell_{G}(E) & = \dim \cF_{G}(E) = \dim (\cR_{G}(E)/
    \ann_{\cR_{G}(E)} \cF_{G}(E))\\
    & = \dim \left (\cR_{G}(E)/\sqrt{\ann_{\cR_{G}(E)}
    \cF_{G}(E)} \right),
    \end{align*}
    the result follows.
\end{proof}

As in the case of ideals, by using the asymptotic value of $\dpt
G_n/E_n$ one can obtain a slightly better bound than the original
one in Burch's inequality. Namely,

\begin{Theorem}\label{bii}
    Let $E \varsubsetneq G \simeq R^e$, $e>0$, be an \rmo{}. Assume $d>0$.
    Then $$ \ell_{G}(E) \leq d+e-1- \dpt(G,E)
    \leq d+e-1- \inf_{n \geq 1} \dpt G_n /E_n . $$
\end{Theorem}

\begin{proof}
    We use induction on $\beta = \dpt (G,E)$ to prove the first
    inequality. If $\beta =0$, we apply
    Corollary \ref{mpG3}. Now, suppose that $\beta> 0$. Let $a
    \in \frm$ such that $a \not \in \z_{R}(G_n/E_n)$ for all
    $n \gg 0$ and $a \not \in \frp$ for all $\frp \in \Min R$, which
    exists by Remark \ref{bi2}.
    By the lemma above
    $$ \ell_G(E) = \ell_{\ol{G}}(\ol{E}), $$
    where $\ol{E}$, $\ol{G}$ are as in Lemma \ref{ad1}.
    Moreover, since $a \not \in \z_{R}(G_n/E_n)$ for all
    $n \gg 0$,
    $$ \dpt G_n/E_n = \dpt (G_n/E_n)/a (G_n/E_n)
    + 1 = \dpt \ol{G}_n/ \ol{E}_n + 1, $$
    for $n \gg 0$. Hence
    $$ \dpt (\ol{G}, \ol{E}) = \dpt (G,E)-1. $$
    Moreover, $\ol{G} \simeq (R/(a))^e$.
    Further, since $a \not \in \frp$ for all $\frp
    \in \Min R$
    $$ \dim R/(a) = \dim R -1 = d-1. $$
    Now by induction,
    $$ \ell_{\ol{G}}(\ol{E}) \leq (d-1) +e-1 - \dpt (\ol{G},
    \ol{E}).$$
    It follows that
    $$ \ell_G(E) = \ell_{\ol{G}}(\ol{E}) \leq d+e-1-
    \dpt (G,E)$$
    as asserted. Finally, since
    $$ \dpt (G,E) \geq \inf_{n \geq 1} \dpt G_n /E_n $$
    the result follows.
\end{proof}

In the case of ideals, it is easy to prove that if if $a \not \in
\bigcup_{n \geq 1} \z_{R}(R/I^n)$ then
$$ \cR(I/aI) \simeq \cR((I+aR)/aR) \simeq \RI \oxr R/aR \simeq 
\RI/a \RI.$$
Furthermore, if $a \in \frm \setminus \bigcup_{n \geq 1}
\z_{R}(R/I^n)$ then
$$  \cF(I/aI) \simeq \RI \oxr R/aR \oxr
    R/\frm \simeq \RI \oxr R/\frm = \FI.$$
For modules we can deduce the following.

\begin{Proposition}\label{bii1}
    Let $E \varsubsetneq G \simeq R^e$, $e>0$, be an \rmo{}.
    Assume $d>0$. If $a \in \frm \setminus \bigcup_{n \geq 1}
    \z_{R}(G_{n}/E_{n})$ then
     \begin{enumerate}
    \item  $ aG_n \cap E_n = a E_n \sube
    \frm E_n$ for all $n \geq 1$;

    \item  $\cR_{\ol{G}}(\ol{E}) \simeq \cR_{G}(E)
    \oxr R/aR \simeq \cR_{G}(E)/a \cR_{G}(E);$

    \item  $\cF_{\ol{G}}(\ol{E}) \simeq \cF_{G}(E).$
    \end{enumerate}
\end{Proposition}

\begin{proof}
    Let $a \in \frm \setminus \bigcup_{n \geq 1}
    \z_{R}(G_{n}/E_{n})$. Hence, as in the
    proof of Lemma \ref{bi3}, we have $aG_n \cap E_n =
    a E_n \sube \frm E_n$ for all $n \geq 1$.
    In particular, $ aG \cap E = a E$. It follows that
    $$ E/aE = E/(aG \cap E) \simeq (E+aG)/aG \simeq
    \ol{E} .$$
    For the assertions (2) and (3), we use Lemma \ref{ad1}. In fact,
    for all $n \geq 1$
    $$\ol{E}_n  =
    (E_n + aG_n)/a G_n \simeq E_n/(aG_n
    \cap E_n) \simeq E_n/a E_n = [\cR_G(E)/a \cR_G(E)]_{n}$$
    and, since $a \in \frm$,
    \begin{align*}
        [\cF_{\ol{G}}(\ol{E})]_{n} & = [\cF_{\ol{G}}(\ol{E})]_{n}
        \simeq E_n / (\frm E_n + (aG_n \cap E_n)) \\
        & \simeq E_n/(\frm E_n + aE_n)
        = E_n/\frm E_n = [\cF_G(E)]_{n}.
    \end{align*}
    The result follows.
 \end{proof}

\section{The case $\cR _G(E)$ being Cohen-Macaulay}

Next we shall prove that, in the case where the Rees algebra of a
\fgtfrmo having rank is \CM, the Burch's inequality is in fact an
equality.

Given a \fg{} module $E$ over a \nor{} $R$ and an $R$-ideal $I$
such that $IE \neq E$
\begin{equation} \label{lch}
    \depth_{I}E = \inf \{ i \in \N_{0} | \h^{i}_{I}(E) \neq 0 \},
\end{equation}
where $\h^{i}_{I}(E)$ denotes the $i$-th local cohomology module
of $E$ \wrt $I$ (cf. \cite[Theorem~6.2.7]{BrSh}). Moreover, if
$\apl{\var}{R}{S}$ is a homomorphism of Noetherian (graded) rings
and $I$ is a (homogeneous) ideal of $R$ and $E$ is a \fg{}
$S$-module
\begin{equation}\label{lch1}
    \h^{i}_{I}(E) \simeq \h^{i}_{IS}(E),
\end{equation}
(cf. \cite[Corollary 35.20]{HeIkOr}).

\begin{Lemma}\label{lc1}
    Let $E \varsubsetneq G \simeq R^e$, $e>0$, be an \rmo{}. Then
    $$ \grd \frm \cR_{G}(E) = \inf_{n \geq 0}  \dpt E_n .$$
\end{Lemma}

\begin{proof}
    We have, by (\ref{lch})
    $$ \grd \frm \cR_G(E) = \depth_{\frm \cR_G(E)} (\cR_G(E))
    = \inf \{ i \in \N_{0} | \h^{i}_{\frm \cR_G(E)} (\cR_G(E))
    \neq 0 \} $$
    since $\cR_G(E)$ is a \nor.
    Moreover, since $\cR_G(E)$ is an \rmo, by (\ref{lch1})
    $$\h^{i}_{\frm \cR_G(E)} (\cR_G(E)) = \h^{i}_{\frm} (\cR_G(E))
    = \h^{i}_{\frm} (\textstyle{\bigoplus_{n \geq 0} E_n})
    = \displaystyle{\bigoplus_{n \geq 0} \h^{i}_{\frm} (E_n)}$$
    by \cite[Theorem~3.4.10]{BrSh}. Therefore
    $$ \h^{i}_{\frm \cR_G(E)} (\cR_G(E)) \neq 0 \iff \exists
     \, m \geq 0 \, \colon \h^{i}_{\frm}(E_{m}) \neq 0. $$
    Now, suppose that $\grd \frm \cR_G(E) = j$.
    Then there exists
    an $m \geq 0$ such that $\h^{j}_{\frm}(E_{m}) \neq 0$, and
    so
    $$ \inf_{n \geq 0}  \dpt E_n  \leq \dpt E_{m}
    \leq j = \grd \frm \cR_G(E). $$
    On the other hand, suppose that $\dpt E_{s} = \inf_{n \geq 0}
     \dpt E_n $ and suppose that \linebreak $\dpt E_{s}
    =j$.
    Then $\h^{j}_{\frm}(E_{s}) \neq 0$, and so
    $\h^{j}_{\frm \cR_G(E)}(\cR_G(E)) \neq 0$. It follows that
    $$ \grd \frm \cR_G(E) \leq j = \dpt E_{s}
    = \inf_{n \geq 0}  \dpt E_n,  $$
    and the equality follows.
\end{proof}

We have the following bound for the depth of the Rees algebra of a
module.

\begin{Proposition}\label{lc2}
    Let $E \varsubsetneq G \simeq R^e$, $e>0$, be an \rmo{}. Then
    $$ \dpt \cR_G(E) \leq \inf_{n \geq 0} \dpt E_n
     + \ell_G(E).$$
\end{Proposition}

\begin{proof}
    By the lemma above we have
    \begin{align*}
        \inf_{n \geq 0}  \dpt E_n  & = \grd \frm \cR_G(E) \\
    & \geq \dpt \cR_G(E) - \dim \cR_G(E)/ \frm \cR_G(E) \\
        & = \dpt \cR_G(E) - \ell_G(E)
    \end{align*}
    by \cite[Th.~17.1]{mat}, proving the inequality.
\end{proof}

The next result was originally proved by D. Eisenbud and C. Huneke
\cite{EiHu} in the ideal case.

\begin{Corollary}\label{lc3}
    Assume that $R$ is Cohen-Macaulay and $d>0$. Let $E\subset G \simeq R^e$,
    $e>0$, be an \rmo{} with rank but not free. If $\cR_G(E)$ is \CM{} then
    $$ \ell_G(E) = d+e-1- \inf_{n \geq 1} \dpt G_n/E_n .$$
\end{Corollary}

\begin{proof}
    By Theorem \ref{bii}
    $$ \ell_G(E) \leq d + e - 1 - \inf_{n \geq 1} \dpt G_n/E_n . $$
    By the proposition above, and since $R$ and $\cR_G(E)$ are
    both Cohen-Macaulay, we have
    \begin{align*}
        \ell_G(E)  & \geq \dpt \cR_G(E) - \inf_{n \geq 0} \dpt E_n \\
    & = \dim \cR_G(E) - (\inf_{n \geq 0} \dpt G_n/E_n + 1) \\
    & = d + e - 1 - \inf_{n \geq 0} \dpt G_n/E_n.
    \end{align*}
    The equality follows.
\end{proof}

Our final result in this section may be useful for induction
arguments.

\begin{Corollary}\label{bii11}
Assume that $R$ is Cohen-Macaulay and $d>0$. Let $E\subset G
\simeq R^e$, $e>0$, be an \rmo{} with rank but not free. If
$\cR_G(E)$ is \CM{} and $\ell_G(E) < d+e-1$ then there exists
     $a \in \frm$ such that
     \begin{enumerate}
    \item  $ aG_n \cap E_n = a E_n \sube
    \frm E_n$ for all $n \geq 1$;
    \item  $\cR_{\ol{G}}(\ol{E}) \simeq \cR_{G}(E)
    \oxr R/aR \simeq \cR_{G}(E)/a \cR_{G}(E);$
    \item  $\cF_{\ol{G}}(\ol{E}) \simeq \cF_{G}(E).$
    \end{enumerate}
\end{Corollary}

\begin{proof} Since $\ell_G(E) < d+e-1$ we have by Corollary
\ref{lc3} that $\inf_{n \geq 1} \dpt G_n/E_n > 0$. Therefore,
$\frm \notin A(n) = \Ass (G_n/E_n)$ for every $n\geq 1$, and so
there exists $a \in \frm \setminus \bigcup_{n \geq 1}
\z_{R}(G_{n}/E_{n})$. Now apply Corollary \ref{bii1}.
 \end{proof}

\section{A criteria for complete intersection}

Let $E\varsubsetneq G \simeq R^e$, $e>0$, be an $R$-submodule of
$G$, with rank $e$ but not free. Following A. Simis, B. Ulrich,
and W. V. Vasconcelos \cite{SiUlVa} we say that $E$ is an ideal
module if the double dual $E^{**}$ is free. Ideal modules provide
a natural extension of several notions in analogy to the case of
ideals. Namely, if $E$ is an ideal module we define the deviation
of $E$ by $\de (E) = \mu (E) -e+1-\Ht F_e(E)$ and the analytic
deviation of $E$ by $\ad (E)=\ell_G(E)-e+1-\Ht F_e(E)$, where $\mu
(\,\cdot\,)$ denotes de minimal number of generators and $F_e(E)$
is the $e$-th Fitting invariant of $E$. Similarly to the ideal
case, one has that the inequalities $\de (E)\geq \ad (E)\geq 0$
hold for any ideal module $E$, see \cite[Proposition 4.2.1]{Co}.
We then say that an ideal module $E$ is a complete intesection if
$\de (E)=0$ and equimultiple if $\ad (E)=0$. Obviously, complete
intersection ideal modules are equimultiple. We also say that an
ideal module $E$ is generically a complete intersection if $\mu
(E_{\frp})=\Ht F_e(E)+e-1$ for all minimal prime ideals $\frp \in
\Min R/F_e(E)$.

\begin{Remark}\label{ci0}
    {\rm Our definitions of deviation and analytic deviation slightly
    differ from those in \cite{SiUlVa} since there it is used
    $\grade F_e(E)$ instead of $\Ht F_e(E)$. Of course, they coincide
    if $R$ is Cohen-Macaulay.}
\end{Remark}

Our aim in this section is to extend to ideal modules the famous
criteria by R. C. Cowsik and M. V. Nori \cite{CoNo} for an ideal
to be a complete intersection. First, and as a consequence of the
Burch's inequality, we have the following criteria for an ideal
module $E$ to be equimultiple in terms of the behaviour of depths
of its Rees powers.

\begin{Corollary}\label{ci1}
    Let $E \varsubsetneq G \simeq R^e$, $e>0$, be an an ideal
    module. If $\depth G_n/E_n=d-\Ht F_e(E)$ for infinitely many
    $n$, then $E$ is equimultiple.
   \end{Corollary}
\begin{proof}
By assumption, $\depth (G,E)=d-\Ht F_e(E)$. Applying Burch's
inequality, we obtain $$\ell _G(E)\leq d+e-1-\depth (G,E)=\Ht
F_e(E)+e-1\leq \ell _G(E)$$ proving that $E$ is equimultiple.
\end{proof}

The following lemma extends to ideal modules a result which has
been very fruitful in the ideal case. We refer to
\cite[Proposition 4.2.14]{Co} for its proof. It follows the same
lines as in the case of ideals by using the theory of reductions
of modules.

\begin{Lemma}\label{ci2}
    Assume that $R$ is Cohen-Macaulay and let $E$ be an ideal
    module. Suppose that $E$ is generically a complete
    intersection. Then $E$ is a complete intersection if and only
    if $E$ is equimultiple.
\end{Lemma}

Now we are ready to prove for ideal modules the extension of the
criteria by R. C. Cowsik and M. V. Nori. We obtain it as
consequence of our version for modules of the Burch's inequality.
In this way, we even get a slightly better version of the
criteria, in the same way as M. Brodmann \cite{Br2} did for
ideals.

\begin{Theorem}\label{ci3}
    Let $R$ be a cohen-Macaulay local ring, $\dim R=d>0$, and let
    $E\subsetneq G\simeq R^e$ be an ideal module having rank $e>0$.
    If $E$ is generically a complete intersecttion then the following
    are all equivalent:

    \begin{enumerate}
    \item[(1)] $E$ is a complete intersection;

    \item[(2)] $G_n/E_n$ are Cohen-Macaulay for all $n>0$;

    \item[(3)] $G_n/E_n$ are Cohen-Macaulay for inifinetly many
    $n$.
    \end{enumerate}

\end{Theorem}

\begin{proof} $(1)\Rightarrow (2)$ This was already proved by
D. Katz and C. Naude \cite[Proposition 3.3]{KaNa} (in fact, they
proved a stronger result: That $G_n/E_n$ are perfect of dimension
$d-\Ht F_e(E)$ for all $n\geq 1$).

$(2)\Rightarrow (3)$ is immediate.

$(2)\Rightarrow (3)$ In virtue of Lemma \ref{ci2} it is enough to
show that $E$ is equimultiple. Since $E$ is generically a complete
intersection, we have by the result of D. Katz and C. Naude that
$\dim G_n/E_n=d-\Ht F_e(E)$ for all $n\geq 1$. Now, by assumption
$$\depth G_n/E_n=\dim G_n/E_n=d-\Ht F_e(E)$$ for infinitely many
$n$, hence by Corollary \ref{ci1} $E$ is equimultiple.
\end{proof}


\end{document}